\documentclass[12pt,a4paper]{amsart}

\usepackage{amsfonts, amsmath, amssymb, amsthm, hyperref}

\usepackage{a4wide}

\newtheorem{thm}{Theorem}
\newtheorem{lem}{Lemma}

\theoremstyle{definition}

\theoremstyle{remark}

\renewcommand{\int}{\mathop{\rm int}}

\renewcommand{\epsilon}{\varepsilon}

\begin{document}
\title{Equipartition of a measure by $(Z_p)^k$-invariant fans}
\author{R.N.~Karasev}
\address{Russia 141700, Dolgoprudny, Institutsky per. 9, Dept. of Mathematics, R.N.~Karasev}
\email{r\_n\_karasev@mail.ru}
\thanks{This research was supported by the Russian Foundation for Basic Research grant 06-01-00648, by the President's of Russian Federation grant MK-1005.2008.1, and partially supported by the Dynasty Foundation.}

\begin{abstract}
We prove a result about partitioning an absolute continuous measure in $\mathbb R^d$ into $2d$ equal parts by a system of cones with common vertex, where $d$ is an odd prime power. The proof is topological and based on the calculation of the equivariant Euler class of a certain vector bundle.
\end{abstract}

\maketitle

\section{Introduction}

Theorems on partitioning a measure or a set of measures in $\mathbb R^d$ into equal parts or parts of prescribed measure have quite a long history, the first to be mentioned is the well-known ``ham sandwich theorem''~\cite{st1942, ste1945}. Some results and problems on measure partitions were proved and mentioned in~\cite{grun1960}. 

A certain type of equipartition results use a partition by hyperplanes or half-hyperplanes, having a common plane of codimension $2$. The results for such partitions have essentially topological proofs, as it can be seen in~\cite{mak1994,barmat2001,zivvre2001}, see also the review~\cite{ziv2004} for different related results. 

Another way of partitioning a measure is to start from a partition of $\mathbb R^d$ by a system of cones with apex at the origin, we call such a system ``fan''. Then this system may be allowed to be moved by translations or by rotations, so as to give a needed partition of a measure, see for example~\cite{zivvre2001}. One case of such results was proved for $\mathbb R^3$ in~\cite{mak1988}, here we generalize it to higher dimensions.

\begin{thm}
\label{mespartcones}
Let $G=(Z_p)^k$, let $d = p^k$ be an odd prime power. Suppose $G$ acts on $\mathbb R^d$ orthogonally by transitively permuting the vectors of some base. Consider a closed cone $C$ with apex at the origin, and suppose that the family of cones $\{\pm g(C)\}_{g\in G}$ gives a partition of $\mathbb R^d$, and the subfamily $\{g(C)\}_{g\in G}$ has a unique common ray.

Then for any absolute continuous probabilistic measure $\mu$ on $\mathbb R^d$ there is a rigid motion $\rho$ that preserves the orientation, and such that for any $g\in G$
$$
\mu(\rho(g(C))) = \mu(\rho(-g(C)) = \frac{1}{2d}.
$$
\end{thm}

For example, we can choose the cones $\pm g(C)$ to be the Voronoi cells of some vector system $\{\pm g(v)\}_{g\in G}$. Note also, that in the general case the cone $C$ need not be convex or polyhedral.

The proof of Theorem~\ref{mespartcones} uses some topological obstruction, that was first used by the author to prove a different result on inscribing a crosspolytope, see Section~\ref{inscrcross} for more comments.

\section{Reduction to a topological fact}

Denote $\{1,2,\ldots,n\} = [n]$.

Consider the group of motions $\mathbb R^d$ that preserve the orientation, denote it $E$. It is natural to write the motions in the from
$$
\rho(x) = \omega (x) + t,
$$
where $\omega\in SO(d)$ is a rotation and $t\in\mathbb R^d$ is the translation vector. Thus, topologically $E=SO(d)\times \mathbb R^d$. 

Now we define a continuous map $\phi$ from $E$ to $V=\mathbb R[G]\oplus\mathbb R[G]$. For $\rho\in E$ and any $g\in G$ put
$$
a_g = \mu(\rho(g(C)),\quad b_g = \mu(\rho(-g(C)),
$$
$(a_g)$ and $(b_g)$ being the coordinates in $V=\mathbb R[G]\oplus\mathbb R[G]$. The measure $\mu$ is absolutely continuous, so the map $\phi$ is continuous. Note that $G$ acts on $E$ by right multiplications by $g^{-1}$, and acts on $V$ by the left multiplication. In these terms the map $g$ is $G$-equivariant (commutes with the $G$-action).

Take new coordinates in $V$ as $s_g = a_g + b_g, t_g = a_g - b_g$. Let us enumerate the elements of $G$ some way as $\{g_1, \ldots, g_d\}$. Consider the one-dimensional subspace $L\subset V$, given by equalities
$$
t_{g_1}=\dots=t_{g_d} = 0,\quad s_{g_1}=\dots=s_{g_d}.
$$
In the space $V/L$ denote the $d$-dimensional linear hull of $\{t_{g_1},\ldots, t_{g_d}\}$ by $U$, and the $d-1$-dimensional linear hull of $\{s_{g_1},\ldots, s_{g_d}\}$ by $W$. Let us denote the natural projection $\pi : V\to V/L$ and $f = \pi\circ \phi$. Now all we need is to prove that the map $f: E\to U\oplus W$ maps some motion $\rho$ to zero. 

\section{Calculating the obstruction}

From here on the cohomology will be considered modulo $p$, its coefficients $Z_p$ are often omitted. For basic facts about (equivariant) topology of vector bundles the reader is referred to the books~\cite{milsta1974,hsiang1975,mishch1998}.

The map $f$ can be regarded as a section of the trivial (but not $G$-trivial) $G$-bundle ever the $G$-space. The first obstruction to the existence of a nonzero section is the Euler class of the bundle. In fact in this problem we could consider the Euler class in the cohomology with compact support, but we use another approach that is more precise.

We can take some $\varepsilon < \frac{1}{2d}$ and assume that for some ball $B'$ centered at the origin $\mu(B') > 1 - \varepsilon$. By the statement of the theorem, the cones $\{\pm g(C)\}$ have only one common point in the origin. Thus, for some big enough ball $B$ the following condition holds: for any $\rho\in SO(d)\times\partial B$ one of the cones $\pm\rho(g(C))$ does not intersect $B'$. Hence, the section $f$ has no zeros on the set $SO(d)\times\partial B \subset E$.

Now we can consider the obstruction as the relative Euler class of $f$ in the cohomology $H_G^{2d-1}(SO(d)\times B, SO(d)\times \partial B)$. In the paper~\cite{kar2008dcpt} the relative Euler class for a section is considered in detail, but here we only need some properties similar to the ordinary Euler class. 

We deform the measure $\mu$ continuously so that it becomes concentrated in $B'$. The section $f$ will be deformed continuously too, remaining nonzero at $SO(d)\times\partial B$, therefore its Euler class does not change.

Let us decompose $f$ into $s_U\oplus s_W$ by the corresponding $G$-bundles. The section $s_U$ has no zeroes on $SO(d)\times \partial B$ and its Euler class is in $e(s_U)\in H_G^d(SO(d)\times B,SO(d)\times \partial B)$, the section $s_W$ has the Euler class $e(s_W)\in H_G^{d-1} (SO(d)\times B) = H_G^{d-1}(SO(d))$. The total Euler class is $e(f) = e(s_U)e(s_W)$.

By the K\"unneth formula 
$$
H_G^*(SO(d)\times B, SO(d)\times \partial B) = H^*(B,\partial B)\otimes H_G^*(SO(d)),
$$
equivalently, $H_G^*(SO(d)\times B, SO(d)\times \partial B) = u \times H_G^*(SO(d))$, where $u$ is the $d$-dimensional generator of $H^*(B,\partial B)$. 

Let us find $e(s_U)$. Fix some rotation $\omega$ in the pair $\rho = (\omega, t)$ and consider the restriction of $s_U$ to $(\{\omega\}\times B, \{\omega\}\times\partial B)$, denote it $s'_U$. Here we formulate a lemma.

\begin{lem}
Under the above assumptions, the map $s'_U$ gives a map of spheres $\partial B\to S(U)$ with nonzero degree mod $p$.
\end{lem}

We postpone the proof of the lemma till the next section. 

It follows from the lemma that the Euler class $e(s'_U)$ is nonzero in $H^*(B,\partial B)$ (the section $s'_U$ cannot be a restriction of some nonzero section over $B$).

In this case $e(s_U)$ under the natural map $H_G^d(B\times SO(d), \partial B\times SO(d))\to H^d(B,\partial B)$ is mapped to nonzero $e(s'_U)\in H^d(B,\partial B)$, this is possible only in the case $e(s_U) = u\times a\in H_G^*(B\times SO(d), \partial B\times SO(d))$, where $a$ is some constant, $a\not=0\mod p$.

Now we use the fact that $e(s_W)$ equals the obstruction that was needed to solve some particular case of Knaster's problem~\cite{kna1947} in the paper~\cite{vol1992}. By its definition $W=\mathbb R[G]/\mathbb R$, this $G$-representation has no irreducible summands with trivial $G$-action, and (see~\cite{hsiang1975}, Chapter III \S 1) its Euler class $e(W)\in H_G^{d-1}(EG) = H^{d-1}(BG)$ is non-zero. Denote the natural classifying $G$-map $p : SO(d)\to EG$. In~\cite{vol1992} (Proposition~3, page~127 in the Russian version) it was shown that the cohomology map $p^* : H_G^k(EG)\to H_G^k(SO(d))$ is injective for every $k < 2(p^k - p^{k-1})$, and it follows that $e(s_W) = p^*e(W)\not=0$.

By the K\"unneth formula and the multiplicative rule for the Euler class we have $e(f)=u\times ae(s_W)\not= 0$. The proof is complete.

\section{Proof of the degree lemma}

Consider the linear subspace $L\in\mathbb R^d$ that corresponds to the unique common ray of the cones $\{g(C)\}_{g\in G}$ and take its orthogonal complement $L^\perp$. $L^\perp$ is a $G$-representation without nonzero fixed points. We can also consider the decomposition $U=U^G\oplus U^\perp$ into the $G$-fixed linear subspace and its orthogonal complement. 

We consider the sphere $S=\partial B$ in $\mathbb R^d$ and spheres $SL = S\cap L, SL^\perp = S\cap L^\perp$. The map $s'_U$ does not map any point in $S$ to zero, so topologically $s'_U$ can be considered as a map to the unit sphere $SU=SU^G*SU^\perp$. The image of the sphere $SL^\perp$ under $s'_U$ does not touch $U^G$, so $s'_U$ induces a $G$-map of spheres $\sigma : SL^\perp\to SU^\perp$. We are going to show that $\sigma$ which has nonzero degree $\mod p$. 

Let us deform the measure $\mu$ to the measure $\mu_0$ that is uniformly distributed in a small ball $B_\varepsilon$ centered at the origin. Under this deformation the map $\sigma$ is deformed into a map $\tau$ of the same degree, if $\varepsilon$ is small enough. Now consider an action of $G$ on $SO(d)\times B$, different from the original action, defined it by $h(\rho) = h^{-1}\rho h$. This is a right action, and it is easy to see that the ``rotation part'' of $\rho$ does not change under this action. Note also that $\mu_0$ is $G$-invariant and
$$
\mu_0(h^{-1}\rho h (\pm g(C))) = \mu_0(\rho h (\pm g(C))),
$$
therefore the original map $s$ and its restriction $\tau : SL^\perp\to SU^\perp$ are also equivariant with respect to this $G$-action (this is true only for the ``symmetric'' measure $\mu_0$). Here it is important that the new action of $G$ takes $\{\omega\}\times SL^\perp$ to itself, unlike the original one. Now we can apply Lemma~2.1 from~\cite{vol2005} that guarantees a nonzero degree mod $p$ of a $G$-equivariant map between spheres of equal dimension.

Now consider some vector $e\in SL$, it is mapped to some unit vector $f\in SU^G$. It can be easily seen that all preimages of $-f$ must be in some small neighborhood of $-e$ (if we chose the ball $B$ large enough). Now consider the half-sphere $H\subset S$ that contains $e$ and has boundary $SL^\perp$. The pair $(H, SL^\perp)$ is mapped to the pair $(SU\setminus\{-f\}, SU\setminus U^G)$, and considering the map between cohomology exact sequences of the pairs, we see that the map
$$
(s'_U)^* : H^d(SU\setminus\{-f\}, SU\setminus U^G, Z_p)=H^d(SU, \{-f\}, Z_p)\to H^d(H, SL^\perp, Z_p) = H^d(S, \{-e\}, Z_p)
$$
is nonzero. This is exactly what we need.

\section{The theorem on inscribing a regular crosspolytope}
\label{inscrcross}

In the paper~\cite{kar2008ins} some results on inscribing a crosspolytope into convex bodies were proved, here we state one of them.

\begin{thm}
\label{insrccross-p}
Let $K\subset\mathbb R^d$ be a smooth convex body, let $d=p^k$ be an odd prime power. Suppose $C$ is a crosspolytope, i.e. the convex hull of vectors $(\pm e_1,\ldots, \pm e_d)$, where $(e_1,\ldots, e_d)$ is some base in $\mathbb R^d$. Suppose the group $G=(Z_p)^k$ acts orthogonally on $\mathbb R^d$ and transitively on the vectors $(e_1,\ldots, e_d)$.

Then there exists a polytope $C'\subset\mathbb R^d$ similar (orientation preserving) to $C$, all its vertices being on $\partial K$.
\end{thm}

In particular, this theorem is true for a regular crosspolytope (formed by an orthonormal base $e_1, \ldots, e_d$) in $\mathbb R^d$, since in this case all the permutations of the base $(e_1, \ldots, e_d)$ are orthogonal transforms of $\mathbb R^d$, and the group of permutations $\mathfrak S_d$ obviously contains a subgroup $(Z_p)^k$.

The proof of Theorem~\ref{insrccross-p} is going to be published in Russian, so it will be useful to show here, how the above reasoning proves this theorem too.

Let us index the vectors $e_i$ by the elements of $G$ as $e_g = ge_1$. We have to consider the configuration space $SO(d)\times K$ and the maps
$$
a_g(p, s) = \max \{a : p + as(e_g)\in K\},\quad b_g = \max \{a : p - as(e_g)\in K\},
$$
they are continuous in the case of smooth strictly convex $K$. Then we again transform them to obtain maps to the same space $U\oplus W$ and note that this map is nonzero on $SO(d)\times \partial K$. The Euler class $e(s_W)$ is the same as above, the Euler class $e(s_U)$ is shown to be nonzero in the similar fashion.

Then we remove the condition of strict convexity by approximating $K$ by strict convex smooth bodies, going to the limit, and showing that the inscribed crosspolytopes do not get degenerate. In~\cite{kar2008ins} it is also shown that $\partial K$ can be replaced by any image of a smooth embedding of a sphere $S^{d-1}\to\mathbb R^d$.

\section{Acknowledgments}

The author thanks the referees for a number of useful remarks that helped to make the presentation of these results much clearer.

\end{document}